\documentclass[12pt,a4paper,oneside,reqno,notitlepage]{amsart}
\usepackage{amssymb}
\textwidth
150mm

\newtheorem{theorem}{Theorem}
\theoremstyle{plain}

\newtheorem{corollary}{Corollary}

\newtheorem{definition}{Definition}

\newtheorem{lemma}{Lemma}

\newtheorem{remark}{Remark}

\numberwithin{equation}{section}

\begin{document}

\baselineskip 8mm
\parindent 9mm

\title[]
{A characteristic of local existence for fractional heat equations in Lebesgue spaces}

\author{ Kexue Li}

\address{Department of Mathematics,
Xi'an Jiaotong University,
 Xi'an
710049, China; Department of Mathematics, University of Washington, Seattle, WA, 98195, USA}
\email{kexueli@gmail.com}

\thanks{{\it 2010 Mathematics Subjects Classification}: 35K55}
\keywords{Fractional heat equation; Dirichlet fractional heat kernel; Local existence; Non-existence.}

\begin{abstract}
In this paper, we consider the fractional heat equation $u_{t}=\triangle^{\alpha/2}u+f(u)$ with Dirichlet boundary conditions on the ball $B_{R}\subset \mathbb{R}^{d}$, where $\triangle^{\alpha/2}$ is the fractional Laplacian, $f:[0,\infty)\rightarrow [0,\infty)$ is continuous and non-decreasing. We present the characterisations of $f$ to ensure the equation has a local solution in $L^{q}(B_{R})$ provided that the non-negative initial data $u_{0}\in L^{q}(B_{R})$. For $q>1$ and $1<\alpha\leq 2$, we show that the equation has a local solution in $L^{q}(B_{R})$ if and only if $\lim_{s\rightarrow \infty}\sup s^{-(1+\alpha q/d)}f(s)=\infty$; and for $q=1$ and $1<\alpha\leq 2$ if and only if $\int_{1}^{\infty}s^{-(1+\alpha/d)}F(s)ds<\infty$, where $F(s)=\sup_{1\leq t\leq s}f(t)/t$. When $\lim_{s\rightarrow 0}f(s)/s<\infty$, the same characterisations holds for the fractional heat equation on the whole space $\mathbb{R}^{d}$.
\end{abstract}
\maketitle

\section{\textbf{Introduction}}
Let $B_{R}:=B_{R}(0)$ be the open ball in $\mathbb{R}^{d}$ of radius $R$ centred at $0$.
We  consider local existence of solutions of the semilinear fractional heat equation
\begin{equation}\label{sfh}
 \ \left\{\begin{aligned}
&u_{t}=\triangle^{\alpha/2}u+f(u),\\
&u(0)=u_{0}\geq 0,
\end{aligned}\right.
\end{equation}
on a domain $\Omega \subset \mathbb{R}^{d}$ and on $\mathbb{R}^{d}$  with Dirichlet boundary conditions, when $\Omega=B_{R}$, $u_{0}\in L^{q}(B_{R})$, $1\leq q<\infty$. $\triangle^{\alpha/2}:=-(-\Delta)^{\alpha/2}$ denotes the fractional Laplacian defined by the Fourier transform
\begin{align*}
(\mathcal{F}(-\triangle)^{\alpha/2}u)(\xi)=|\xi|^{\alpha}\mathcal{F}(u)(\xi),
\end{align*}
where $0<\alpha\leq2$, $\mathcal{F}$ denotes the Fourier transform. If $\alpha=2$, the fractional Laplacian $\triangle^{\alpha/2}$ becomes the Laplacian $\triangle$.

We present the characteristic conditions on $f$ such that (\ref{sfh}) has a local solution bounded in $L^{q}(\Omega)$ for non-negative initial data in $L^{q}(\Omega)$. Many previous results about the existence of semilinear parabolic equations or semilinear fractional parabolic equations are concerned with the nonlinearity term $f(u)=u^{p}$, or $f$ is assumed to be convex. In this paper, we only require that $f: [0,\infty)\rightarrow [0,\infty)$ is continuous and non-decreasing. The non-existence results are based on lower bounds of the fractional heat kernel, particularly on lower bounds of the function that the fractional heat semigroup acts on the characteristic function of a ball.

For $\alpha=2$,  when $\Omega=\mathbb{R}^{d}$, $f(u)=u^{p}(p>1)$, Fujita \cite{Fujita} proved that (\ref{sfh}) has a global solution $u(t,x)$ for sufficiently small $u_{0}$ if $d>\frac{2}{p-1}$, but (\ref{sfh}) has no global solution for any $u_{0}\not\equiv 0$ if $d>\frac{2}{p-1}$. For the critical case $d=\frac{2}{p-1}$,  Hayakawa \cite{Hayakawa} showed that there is no global solution if $d$ equals to 1 or 2. Kobayashi et al. \cite{KST} obtained the same result for general $d$. The critical exponent for nonlinear evolution equations has been studied by many authors, see for example, \cite{Levine}, \cite{Deng}, and references therein.  Weissler \cite{Weissler} studied the semilinear parabolic equation
\begin{align}
u_{t}=Au+f(u), \ u(0)=u_{0},
\end{align}
where $u: [0,T]\rightarrow L^{q}(\Omega)$, $\Omega$ is a domain in $\mathbb{R}^{d}$ and $1\leq q<\infty$,  $A$ is the infinitesimal generator of a $C_{0}$ semigroup $e^{tA}$ on $L^{q}(\Omega)$, $f$ has a polynomial growth. One main difference between \cite{Weissler} and previous work is that the initial data $u(0)$ is allowed to be an arbitrary function $u_{0}$ in $L^{q}(\Omega)$. If $\alpha=2$, $f(u)=|u|^{p-1}u$, (\ref{sfh}) was considered by Brezis and Cazenave \cite{BC} when $\Omega\subset \mathbb{R}^{d}$ is smooth bounded domain and $p>1$, $u_{0}\in L^{q}(\Omega)$ for $1\leq q<\infty$. The critical exponent is $q=d(p-1)/2$. They obtained the existence and uniqueness of a local solution for any $u_{0}\in L^{q}(\Omega)$ provided that $q>d(p-1)/2$ $(resp. \ q=d(p-1)/2)$ and $q\geq 1$ $(resp. \ q>1)$, $d\geq1$.

In the last decades, fractional laplacian has attracted a great deal of attention, which naturally arises in anomalous diffusion \cite{MDH}, conformal geometry \cite{Chang},  phase transitions \cite{Sire}, quasi-geostrophic flows \cite{CV}, etc. Differential equations with fractional Laplacian appeared in many papers, see for example \cite{BCDS}, \cite{RE}, \cite{CDDF}, \cite{Fino}, \cite{Silvestre}.
Sugitani \cite{Sugitani} studied the equation
\begin{equation}\label{sugitani}
 \ \left\{\begin{aligned}
&u_{t}=-(-\frac{\triangle}{2})^{\alpha/2}u+f(u),\\
&u(0)=u_{0},
\end{aligned}\right.
\end{equation}
where $0<\alpha\leq 2$, $u_{0}$ is nontrivial, nonnegative and continuous function on $\mathbb{R}^{d}$. Let $f(u)$ be a nonnegative continuous function with $f(0)=0$, defined on $[0,\infty)$, satisfying the following conditions: \\
(f.1) $f$ is increasing and convex. \\
(f.2) There exists some $\beta\in [0, \frac{\beta}{d}]$ and $c'\in (0,\infty)$, such that $\lim_{u\downarrow 0}\frac{f(u)}{u^{1+\beta}}=c'$.\\
(f.3) $\int_{1}^{\infty}\frac{du}{f(u)}<\infty$. \\
Under the above assumptions, Sugitani proved that the nonnegative solution $u(t,x)$ of (\ref{sugitani}) blows up in finite time, i.e., there exists some $t_{0}>0$ such that $u(t,x)=\infty$ for every $t\geq t_{0}$ and $x\in \mathbb{R}^{d}$. Fino and Kirane \cite{FK} considered a time-space nonlocal nonlinear parabolic equation
\begin{equation}\label{FK}
 \ \left\{\begin{aligned}
&u_{t}+(-\triangle)^{\alpha/2}u=\frac{1}{\Gamma(1-\beta)}\int_{0}^{t}(t-s)^{-\beta}|u|^{p-1}u(s)ds, \ x\in \mathbb{R}^{d}, \ t>0,\\
&u(0)=u_{0},
\end{aligned}\right.
\end{equation}
where $d\geq 1$, $0<\alpha\leq 2$, $0<\beta<1$, $p>1$ and $u_{0}\in C_{0}(\mathbb{R}^{d})$, which denotes the space of all continuous functions tending to zero at infinity. They showed that blowing-up solution exist and study their time blowup profile. The necessary conditions for local or global existence is established. Tan and Xu \cite{ZY} considered the problem
\begin{equation}\label{zy}
 \ \left\{\begin{aligned}
&u_{t}=-(-\triangle)^{\alpha}u+h(t)u^{v}, \ (x,t)\in \mathbb{R}^{d}\times (0,T),\\
&u(0)=u_{0}, \ u_{0}\geq 0, \ u_{0}\not\equiv 0,
\end{aligned}\right.
\end{equation}
where $0<\alpha\leq 1$, $v>1$, $u$ is a curve in $L^{q}(\mathbb{R}^{d})(1<q<\infty)$, $u: [0,T]\rightarrow L^{q}(\mathbb{R}^{d})$, $h(t)$ is continuous and $c_{0}t^{\sigma}\leq h(t)\leq c_{1}t^{\sigma}$ for large $t$, $c_{0},\ c_{1}>0$, $\sigma>-1$ are constants. The initial value $u(0)$ is assumed to be in $L^{q}(\mathbb{R}^{d})$. The authors discussed the existence and non-existence of global solutions.

However, if we only assume that $f$ is monotonic, most of the above existence results about semilinear parabolic equations and semilinear fractional parabolic equations will be not valid. Laister et al. \cite{LRSV} considered the semilinear heat equation $u_{t}-\triangle u=f(u)$, where $f: [0,\infty)\rightarrow [0,\infty)$ is continuous and non-decreasing but need not be convex. For $q\geq 1$, the authors gave the characteristic conditions for $f$ such that the equation has a local solution bounded in $L^{q}(\Omega)$ for all non-negative initial value $u_{0}\in L^{q}(\Omega)$, where $\Omega\subset \mathbb{R}^{d}$ is a bounded domain with Dirichlet boundary conditions.  They also presented the characteristic conditions for the equation on the whole space $\mathbb{R}^{d}$. It is pointed out in \cite{LRSV} that local non-existence is `immediate blowup' in some sense, which is different from `finite-time blowup'(see, \cite{Kaplan}, \cite{Ni}, \cite{LN}).

In this paper, we give a characteristic of $f$ such that for $q\geq 1$, (\ref{sfh}) has a local solution bounded in $L^{q}(B_{R})$ for all non-negative initial value $u_{0}\in L^{q}(B_{R})$. We show that if for $q\geq 1$, $1<\alpha\leq 2$, $\lim_{s\rightarrow \infty}\sup s^{-(1+\alpha q/d)}f(s)=\infty$, then there exists a non-negative initial value $u_{0}$ for which (\ref{sfh}) has no local solution bounded in $L^{q}(B_{R})$. If $\lim_{s\rightarrow \infty}\sup s^{-(1+\alpha q/d)}f(s)<\infty$, we see that $f(s)\leq c(1+s^{1+\alpha q/d})$ for some positive constant $c$. From the fact that $u_{t}=\triangle^{\alpha/2}u+c(1+u^{p})$, $p=1+\alpha q/d$ has the local $L^{q}$ existence property, it follows (see Theorem \ref{property}) that for $q>1$ and $1<\alpha\leq 2$, equation (\ref{sfh}) has the local $L^{q}$ existence property for every non-negative $u_{0}\in L^{q}(B_{R})$ if and only if $\lim_{s\rightarrow \infty}\sup s^{-(1+\alpha q/d)}f(s)<\infty$. The case $q=1$ is more delicate. We prove that the condition $\sum_{k=1}^{\infty}s_{k}^{-(1+\alpha/d)}f(s_{k})=\infty$  for some sequence $s_{k}$ and $1<\alpha\leq2$ such that $s_{k+1}\geq \tau s_{k} \ (\tau>1)$ can guarantee a non-existence result. By the equivalent integral condition $\int_{1}^{\infty}s^{-(1+\alpha/d)}F(s)=\infty$, $F(s)=\sup_{1\leq t\leq s}\frac{f(t)}{t}$ and some properties of supersolutions, we obtain a result that equation (\ref{sfh}) has the local $L^{1}$ existence property if and only if
\begin{align*}
\int_{1}^{\infty}s^{-(1+\alpha/d)}F(s)ds<\infty, \ \mbox{where} \ F(s)=\sup_{1\leq t\leq s}\frac{f(t)}{t}.
\end{align*}

The paper is organized as follows. In Section 2, we show some lower bounds on solutions of fractional heat equation for the initial value is the characteristic function of a ball.
In Section 3, we consider the problem with the initial data in $L^{q}(B_{R})$, $1<q<\infty$. Section 4 contains the results for initial data in $L^{1}(B_{R})$. In Section 5, we discuss the problem on the whole space $\mathbb{R}^{d}$.

\section{Lower bounds on solutions of the Dirichlet fractional heat equation}
For any $r>0$, we denote by $\chi_{r}$  the characteristic function of $B_{r}:=B_{r}(0)$. The volume of the unit ball in $\mathbb{R}^{d}$ is $\omega_{d}$.
The solution of the fractional heat equation on $B_{R}$ with Dirichlet boundary condition
\begin{equation}\label{fhe}
 \ \left\{\begin{aligned}
&u_{t}=\triangle^{\alpha/2}u, \\
& u(x,0)=u_{0}(x)\in L^{1}(B_{R}), \\
&u|_{\partial B_{R}}=0,
\end{aligned}\right.
\end{equation}
can be expressed as
\begin{equation}\label{solution}
u(t,x)=(S_{\alpha}(t)u_{0})(x):=\int_{B_{R}}p_{D}(t,x,y)u_{0}(y)dy,
\end{equation}
where $p_{D}(x,y,t)$ is the heat kernel of $\Delta^{\alpha/2}$ on $B_{R}$ with Dirichlet boundary $u|_{\partial B_{R}}=0$. It is well known that  $S_{\alpha}(t)$ is a semigroup with the generator $\Delta^{\alpha/2}$, see \cite{KT}.

Let $p(t,x,y)$ be the heat kernel of $\Delta ^{\alpha/2}$ on $\mathbb{R}^{d}$.
We have the following inequality (see \cite{BGR})
\begin{align}\label{inequality}
0\leq p_{D}(t,x,y)\leq p(t,x,y) \ {\mbox{for all}} \ t>0, \ x,y \in \mathbb{R}^{d}.
\end{align}
Throughout this paper, we use $c_{0}$, $c_{1}$, $c_{2}$, $\ldots$ to denote generic constants, which may change from line to line.
For two nonnegative functions $f_{1}$ and $f_{2}$, the notion $f_{1}\asymp f_{2}$ means that $c_{1}f_{2}(x)\leq f_{1}(x)\leq c_{2}f_{2}(x)$, where $c_{1},c_{2}$ are positive constants.
It is well known that (see, e.g., \cite{TG,KT,CKS})
\begin{align*}
p(t,x,y)\asymp \big(t^{-d/\alpha}\wedge \frac{t}{|x-y|^{d+\alpha}}\big),
\end{align*}
that is, there exist constants $c_{1}$, $c_{2}$ such that for $t>0$, $x,y\in \mathbb{R}^{d}$,
\begin{align}\label{c1}
c_{1}\big(t^{-d/\alpha}\wedge \frac{t}{|x-y|^{d+\alpha}}\big)\leq p(t,x,y)\leq c_{2}\big(t^{-d/\alpha}\wedge \frac{t}{|x-y|^{d+\alpha}}\big),
\end{align}
where $c_{1}$ and $c_{2}$ are positive constants depending on $\alpha$. \\
From (\ref{c1}), it is easy to see that $p(t,x,y)$ satisfies the following inequality
\begin{align}\label{alpha}
\frac{c_{1}t}{(t^{1/\alpha}+|y-x|)^{d+\alpha}}\leq p(t,x,y)\leq \frac{c_{2}t}{(t^{1/\alpha}+|y-x|)^{d+\alpha}},
\end{align}
where $t>0$, $x,y\in \mathbb{R}^{d}$, $c_{1}$ and $c_{2}$ are positive constants depending on $\alpha$.
\begin{lemma}\label{Dirichlet}
(see \cite{Chen}, Theorem 1.1.) Let $D$ be a $C^{1,1}$ open subset of $\mathbb{R}^{d}$ with $d\geq 1$ and $\delta_{D}(x)$ the Euclidean distance between $x$ and $D^{c}$. \\
(i) For every $T>0$, on $(0,T]\times D\times D$,
\begin{align*}
p_{D}(t,x,y)\asymp \big(1\wedge \frac{\delta_{D}(x)^{\alpha/2}}{\sqrt{t}}\big)\big(1\wedge \frac{\delta_{D}(y)^{\alpha/2}}{\sqrt{t}}\big)\big(t^{-1/\alpha}\wedge \frac{t}{|x-y|^{1+\alpha}}\big).
\end{align*}
(ii) Suppose in addition that $D$ is bounded. For every $T>0$, there are positive constants $c_{1}< c_{2}$ such that on $[T,\infty)\times D\times D$,
\begin{align*}
c_{1}e^{-\lambda_{0}t}\delta_{D}(x)^{\alpha/2}\delta_{D}(y)^{\alpha/2}\leq p_{D}(t,x,y)\leq c_{2}e^{-\lambda_{0}t}\delta_{D}(x)^{\alpha/2}\delta_{D}(y)^{\alpha/2},
\end{align*}
where $\lambda_{0}>0$ is the smallest eigenvalue of the Dirichlet fractional Laplacian $(-\Delta)^{\alpha/2}|_{D}$, $p_{D}(t,x,y)$ is the Dirichlet fractional heat kernel.
\end{lemma}
For the definition of $C^{1,1}$ open set, we refer to  \cite{Chen}.

\begin{lemma}\label{constant}
For any $r>0$, $\delta>0$ for which $B_{r+2\delta}\subset B_{R}$, there exists a positive constant $c(d,\alpha)$, such that for all $0<t\leq \delta^{\alpha}$,
\begin{align}\label{dalpha}
S_{\alpha}(t)\chi_{r}\geq c(d,\alpha)\big(\frac{r}{r+t^{1/\alpha}}\big)^{d}\chi_{r+t^{1/\alpha}}.
\end{align}
\end{lemma}
\begin{proof}
Since $B_{r+2\delta}\subset B_{R}$, for all $t\in (0,\delta^{\alpha})$ and $x\in B_{r+t^{1/\alpha}}$, we have $\mbox{dist}(x, B(0,R)^{c})\geq \delta$.
For such $x$ and $t$, by (i) of Lemma \ref{Dirichlet},
\begin{align}\label{Chen}
p_{D}(t,x,y)\geq c_{0}(t^{-d/\alpha}\wedge \frac{t}{|x-y|^{d+\alpha}}).
\end{align}
Then
\begin{align*}
[S_{\alpha}(t)\chi_{r}](x)&=\int_{B_{r}(0)}p_{D}(t,x,y)dy\nonumber\\
&\geq c_{0}\int_{B_{r}}(t^{-d/\alpha}\wedge \frac{t}{|x-y|^{d+\alpha}})dy.
\end{align*}
And we have
\begin{align}\label{chi}
[S_{\alpha}(t)\chi_{r}](x) \geq c_{1}\int_{B_{r}}\frac{t}{(t^{1/\alpha}+|y-x|)^{d+\alpha}}dy.
\end{align}
Since $\int_{B_{r}}\frac{t}{(t^{1/\alpha}+|y-x|)^{d+\alpha}}dy$ is radially symmetric and decreasing with $|x|$, then for $|x| \leq r+t^{1/\alpha}$ and any unit vector $\tau$, we have
\begin{align}\label{solution}
[S_{\alpha}(t)\chi_{r}](x) &\geq c_{1}\int_{B_{r}}\frac{t}{(t^{1/\alpha}+|y-x|)^{d+\alpha}}dy\nonumber\\
&\geq c_{1}\int_{B_{r}((r+t^{1/\alpha})\tau)}\frac{t}{(t^{1/\alpha}+|z|)^{d+\alpha}}dz\nonumber\\
&=c_{1}\int_{B_{rt^{-1/\alpha}}((t^{-1/\alpha}r+1)\tau)}\frac{1}{(1+|w|)^{d+\alpha}}dw.
\end{align}
Note that $B_{t^{-1/\alpha}r}((t^{-1/\alpha}r+1)\tau)\subset B_{t^{-1/\alpha}\rho}((t^{-1/\alpha}\rho+1)\tau)$, if $\rho\geq r$.  Then for $r\geq t^{1/\alpha}$,  by (\ref{solution}), we have
\begin{align}\label{geq}
[S_{\alpha}(t)\chi_{r}](x)&\geq c_{1}\int_{B_{1}(2\tau)}\frac{1}{(1+|w|)^{d+\alpha}}dw\nonumber\\
&=c_{2}(d,\alpha),
\end{align}
where $c_{2}(d,\alpha)$ is a positive constant, which only depends on $d$ and $\alpha$.
For $r\leq  t^{1/\alpha}$, by (\ref{solution}), we get
\begin{align}\label{leq}
[S_{\alpha}(t)\chi_{r}](x)&\geq c_{1}\int_{B_{rt^{-1/\alpha}}((t^{-1/\alpha}r+1)\tau)}\frac{1}{(2+2rt^{-1/\alpha})^{d+\alpha}}dw\nonumber\\
&\geq \frac{c_{1}(rt^{-1/\alpha})^{d}}{(2+2rt^{-1/\alpha})^{d+\alpha}}\nonumber\\
&\geq c_{3}(d,\alpha)(rt^{-1/\alpha})^{d},
\end{align}
where $c_{3}(d,\alpha)$ is a positive constant, which only depends on $d$ and $\alpha$. \\
Set $c(d,\alpha)=\min\{c_{2}(d,\alpha), c_{3}(d,\alpha)\}$. Then by (\ref{geq}) and (\ref{leq}),
\begin{align*}
[S_{\alpha}(t)\chi_{r}](x)&\geq c(d,\alpha)\big(\frac{r}{\max\{r,t^{1/\alpha}\}}\big)^{d}\\
&\geq c(d,\alpha)\big(\frac{r}{r+t^{1/\alpha}}\big)^{d}.
\end{align*}
\end{proof}
\begin{corollary}\label{integrate}
For any $r,\delta>0$ for which $B_{r+2\delta}\subset B_{R}$, there exists a constant $\mu(d,\alpha)$, depending only on $d$ and $\alpha$, such that
for all $0<t\leq \delta^{\alpha}$,
\begin{align*}
\int_{B_{R}}S_{\alpha}(t)\chi_{r}dx\geq \mu(d,\alpha)r^{d}.
\end{align*}
\end{corollary}

\begin{proof}
By (\ref{dalpha}),
\begin{align*}
\int_{B_{R}}S_{\alpha}(t)\chi_{r}dx&\geq c(d,\alpha)\big(\frac{r}{r+t^{1/\alpha}}\big)^{d}\int_{B_{R}}\chi_{r+t^{1/\alpha}}dx\\
&\geq c(d,\alpha)\omega_{d}r^{d}.
\end{align*}
\end{proof}

\begin{corollary}\label{for which}
For any $r,\delta>0$ for which $B_{r+2\delta}\subset B_{R}$, there exists a constant $\nu(d,\alpha)$, depending only on $d$ and $\alpha$, such that
for all $0<t\leq \min\{\delta^{\alpha},r^{\alpha}\}$,
\begin{align*}
S_{\alpha}(t)\chi_{r}\geq \nu(d,\alpha)\chi_{r+t^{1/\alpha}}.
\end{align*}
\end{corollary}

\section{Initial data in $L^{q}(B_{R})$, $1<q<\infty$}
\begin{definition}
Assume that $f:[0,\infty)\rightarrow [0,\infty)$ and $u_{0}\geq 0$, $u$ is said to be a local integral solution of (\ref{sfh}) on $[0,T)$ for some $T>0$ if $u: B_{R}\times [0,T)\rightarrow [0,\infty)$ is measurable, finite almost everywhere, and satisfies
\begin{align}\label{measurable}
u(t)=S_{\alpha}(t)u_{0}+\int_{0}^{t}S_{\alpha}(t-s)f(u(s))ds
\end{align}
almost everywhere in $B_{R}\times [0,T)$.
\end{definition}
\begin{definition}
$u$ is called a local $L^{q}$ solution of (\ref{sfh}) if $u$ is a local integral solution on $[0,T)$ for some $T>0$ and $u\in L^{\infty}((0,T); L^{q}(B_{R}))$. If for every non-negative $u_{0}\in L^{q}(B_{R})$, the integral solution $u\in L^{q}(B_{R})$, we say that (\ref{sfh}) has the local existence property in $ L^{q}(B_{R})$.
\end{definition}
\begin{theorem}\label{local solution}
Let $f: [0,\infty)\rightarrow [0,\infty)$ be non-decreasing. If $q\in [1,\infty)$, $\alpha\in (1,2]$ and
\begin{align}\label{infty}
\lim_{s\rightarrow \infty}\sup s^{-(1+\alpha q/d)}f(s)=\infty,
\end{align}
there exists a non-negative ${u}_{0}\in L^{q}(B_{R})$ such that
\begin{equation}\label{initial}
 \ \left\{\begin{aligned}
&u_{t}=\triangle^{\alpha/2}u+f(u), \\
& u(x,0)=u_{0}(x), \\
&u|_{\partial B_{R}}=0
\end{aligned}\right.
\end{equation}
has no local $L^{q}$ solution.
\end{theorem}
\begin{proof}
 By (\ref{infty}), there exists an increasing sequence $\phi_{k}(k=1,2,\cdots)$ such that
\begin{equation*}
\phi_{k}\geq k, \ \ f(\phi_{k})\geq \phi_{k}^{p}e^{k/q},
\end{equation*}
where $p=1+\alpha q/d$. \\
Set $r_{k}=\varepsilon \phi_{k}^{-q/d}k^{-\alpha q/d}$ and express the initial data $u_{0}(x)=\sum_{k=1}^{\infty}u_{k}$, $u_{k}=\frac{1}{\nu(d,\alpha)}\phi_{k}\chi_{r_{k}}$,
where $\nu(d,\alpha)$ is the same constant as that in Corollary \ref{for which}, $\varepsilon$ is small enough such that $B_{3r_{k}}\subset B_{R}$ for every $k$. Since
\begin{align*}
\|u_{k}\|_{L^{q}}^{q}&=\frac{\omega_{d}}{\nu(d,\alpha)^{q}}r_{k}^{d}\phi_{k}^{q}\\
&=\frac{\omega_{d}}{\nu(d,\alpha)^{q}}(\varepsilon \phi_{k}^{-q/d}k^{-\alpha q/d})^{d}\phi_{k}^{q}\\
&=\frac{\omega_{d}}{\nu(d,\alpha)^{q}}\varepsilon^{d}k^{-\alpha q},
\end{align*}
then
\begin{align*}
\|u_{k}\|_{L^{q}}=\frac{\omega_{d}^{1/q}}{\nu(d,\alpha)}\varepsilon^{d/q}k^{-\alpha},
\end{align*}
and
\begin{align*}
\|u_{0}\|_{L^{q}}\leq \sum_{k=1}^{\infty}\|u_{k}\|_{L^{q}}=\frac{\omega_{d}^{1/q}}{\nu(d,\alpha)}\varepsilon^{d/q}\sum_{k=1}^{\infty}k^{-\alpha}<\infty.
\end{align*}
If a solution $u(t)$ of (\ref{initial}) exists, then
\begin{align}\label{then}
u(t)=S_{\alpha}(t)u_{0}+\int_{0}^{t}S_{\alpha}(t-s)f(u(s))ds.
\end{align}
From $u\geq 0$, $f\geq 0$, it follows that
\begin{align}\label{that}
u(t)\geq S_{\alpha}(t)u_{0}\geq S_{\alpha}(t)u_{k}.
\end{align}
From now, we fix $k$. By (\ref{then}), (\ref{that}) and  note that $f$ is non-decreasing,
\begin{align}\label{lower bound}
u(t)\geq \int_{0}^{t}S_{\alpha}(t-s)f(s_{\alpha}(s)u_{k})ds.
\end{align}
In corollary \ref{for which}, taking $\delta=r_{k}$, we get
\begin{align}\label{alphak}
S_{\alpha}(s)u_{k}&=S_{\alpha}(s)\big(\frac{1}{\nu(d,\alpha)}\phi_{k}\chi_{r_{k}}\big)\nonumber\\
&=\frac{\phi_{k}}{\nu(d,\alpha)}S_{\alpha}(s)\chi_{r_{k}}\nonumber\\
&\geq \phi_{k}\chi_{r_{k}+t^{1/\alpha}}\nonumber\\
&\geq \phi_{k}\chi_{r_{k}},
\end{align}
where $s\in [0,r_{k}^{1/\alpha}]$. \\
By (\ref{alphak}), we have
\begin{align*}
f(S_{\alpha}(s)u_{k})\geq f(\phi_{k})\chi_{r_{k}}, \ \ 0<s\leq r_{k}^{\alpha}.
\end{align*}
This together with Corollary \ref{for which} yield
\begin{align}\label{together}
S_{\alpha}(t-s)f(S_{\alpha}(s)u_{k})&\geq \nu(d,\alpha)f(\phi_{k})\chi_{r_{k}+(t-s)^{1/\alpha}}\nonumber\\
&\geq \nu(d,\alpha)f(\phi_{k})\chi_{r_{k}},
\end{align}
where $0<s\leq t\leq r_{k}^{\alpha}$.\\
Set $t_{k}=r_{k}^{\alpha}$, for any $t\in [t_{k}/2, t_{k}]$, by (\ref{lower bound}), (\ref{together}), we have
\begin{align*}
[u(t)](x)&\geq \int_{0}^{t}S_{\alpha}(t-s)f(S_{\alpha}(s)u_{k})ds\\
&\geq \int_{0}^{r_{k}^{\alpha}/2}\nu(d,\alpha)f(\phi_{k})\chi_{r_{k}}ds\\
&\geq \frac{1}{2}r_{k}^{\alpha}\nu(d,\alpha)f(\phi_{k})\chi_{r_{k}}.
\end{align*}
Therefore
\begin{align*}
\|u(t)\|_{L^{q}}^{q} &\geq \int_{B_{r_{k}}}|u(t)|^{q}dx\\
&\geq \big(\frac{\nu(d,\alpha)}{2}\big)^{q}r_{k}^{\alpha q}f(\phi_{k})^{q}r_{k}^{d}\omega_{d}\\
&=\big(\frac{\nu(d,\alpha)}{2}\big)^{q}\omega_{d}r_{k}^{\alpha q+d}f(\phi_{k})^{q}\\
&\geq \big(\frac{\nu(d,\alpha)}{2}\big)^{q}\omega_{d}(\varepsilon \phi_{k}^{-q/d}k^{-\alpha q/d})^{\alpha q+d}\phi_{k}^{(1+(\alpha q/d))q}e^{k}\\
&\geq \big(\frac{\nu(d,\alpha)}{2}\big)^{q}\omega_{d}\varepsilon^{\alpha q+d}k^{-\alpha q(1+\alpha q/d)}e^{k}\\
&\rightarrow \infty
\end{align*}
as $k\rightarrow \infty$. This implies $u$ does not belong to $L^{\infty}(0,T); L^{q}(B_{R}))$ for any $T>0$.
\end{proof}

The initial data is always assumed to be an element $\mathcal{M}_{S}^{+}(B_{R})$, the set of nonnegative,
a.e. finite measurable functions on $B_{R}$ such that
\begin{align*}
S_{\alpha}(t)u_{0}=\int_{B_{R}}p_{D}(t,\cdot,y)u_{0}(y)dy<\infty \ \  \mbox{for} \ \  t>0.
\end{align*}
Define the operator
\begin{align}\label{equation}
\mathcal{F}[v](t)=S_{\alpha}(t)u_{0}+\int_{0}^{t}S_{\alpha}(t-s)f(v(s))ds.
\end{align}

Let $\mathcal{M}^{+}(B_{R})$ be the set of nonnegative, almost everywhere finite, measurable
functions on $B_{R}$. A function $u$ is said to be a solution of (\ref{equation}) if for any $u\in \mathcal{M}^{+}(B_{R})$ such that
$\mathcal{F}[u]=u$ a.e. in  $B_{R}$. Any $w\in \mathcal{M}^{+}(B_{R})$  satisfying $\mathcal{F}[w]\leq w$ (resp. $\mathcal{F}[w]\geq w$)
will be called a supersolution (resp. subsolution).
\begin{lemma}\label{solution}
Suppose $f: [0,\infty)\rightarrow [0,\infty)$ is continuous, nondecreasing and let $u_{0}\in \mathcal{M}_{S}^{+}(B_{R})$. Then $\mathcal{\mathcal{F}}$ admits a solution in $B_{R}$ if and only if it admits a
supersolution in $B_{R}$.
\end{lemma}
\begin{proof}
The proof  is essentially the same as that of Theorem 1 in \cite{Robinson}, with the difference that $S_{\alpha}(t)$ is changed to $S(t)$, the heat semigroup. So it is omitted.
\end{proof}

\begin{theorem}\label{property}
Let $f:[0,\infty)\rightarrow [0,\infty)$ be non-decreasing and continuous. If $q\in (1,\infty)$, $\alpha\in(1,2]$, then (\ref{initial}) has the local existence property in $L^{q}(B_{R})$ if and only if
\begin{align*}
\lim_{s\rightarrow \infty}\sup s^{-(1+\alpha q/d)}f(s)<\infty.
\end{align*}
\end{theorem}
\begin{proof}
By Theorem \ref{local solution}, we only need to prove the sufficiency. Since $\lim_{s\rightarrow \infty}\sup s^{-(1+\alpha q/d)}f(s)<\infty$, there exists a positive constant $c$ such that
\begin{align}\label{comparion}
f(s)\leq c(1+s^{p}), \ \mbox{where}\ p=1+\alpha q/d.
\end{align}
Consider the equation
\begin{align}\label{analytic}
u_{t}=\triangle^{\alpha/2}u+c(1+u^{p}).
\end{align}
Set $A=-\Delta$, it is known that $-A$ generates a strongly continuous bounded semigroup $S(t)=\exp(-tA)$ on $L^{p}(B_{R})$ (see Theorem 3.5 in \cite{Pazy}). Then the fractional power $-A^{\alpha/2}=\Delta^{\alpha/2}$, $0<\alpha<2$, generates an analytic semigroup $S_{\alpha}(t)=\exp(-tA^{\alpha/2})=\exp(t\Delta^{\alpha/2})$ on $L^{p}(B_{R})$ (See \cite{Yosida}). Then we have the following smoothing effect (See \cite{Fino},  for example)
\begin{align}\label{smooth effect}
\|S_{\alpha}(t)\varphi\|_{L^{q}}\leq Ct^{-\frac{d}{\alpha}(\frac{1}{r}-\frac{1}{q})}\|\varphi\|_{L^{r}},
\end{align}
for all $\varphi\in L^{r}$ and all $1\leq r\leq q\leq \infty, \ t>0$.
Similar to the proof of Theorem 2.2 in \cite{ZY} (or the proof of Theorem 6.1 in \cite{FK}),  we can get the local $L^{q}$ existence property of (\ref{analytic}).  This together with Lemma \ref{solution} yield the conclusion.
\end{proof}

\section{Initial data in $L^{1}(B_{R})$}
\subsection{A condition for non-existence of a local $L^{1}$ solution}
\begin{theorem}\label{assume}
Assume that $f: [0,\infty)\rightarrow [0,\infty)$ is non-deceasing and there exists a sequence $\{s_{k}\}_{k=1}^{\infty}$ such that
\begin{align*}
s_{k+1}\geq \tau s_{k},
\end{align*}
and
\begin{align*}
\sum_{k=1}^{\infty}s_{k}^{-p}f(s_{k})=\infty,
\end{align*}
where $\tau>1$, $p=1+\frac{\alpha}{d}$, $\alpha\in(1,2]$. Then there exists a non-negative initial data $u_{0}\in L^{1}(B_{R})$ such that
\begin{equation}\label{L1}
 \ \left\{\begin{aligned}
&u_{t}=\triangle^{\alpha/2}u+f(u), \\
& u(x,0)=u_{0}(x), \\
&u|_{\partial B_{R}}=0
\end{aligned}\right.
\end{equation}
has no local integral solution that belongs to $L^{1}_{loc}(B_{R})$ for $t>0$, thus no local $L^{1}$ solution exists.
\end{theorem}
\begin{proof}
Let $\phi_{k}=\frac{s_{k}}{c(d,\alpha)}$ and $u_{n}(x)=\frac{1}{n^{\alpha}}\beta_{n}^{d}\chi_{1/\beta_{n}}$, where $n=1,2,\cdots$, $\beta_{n}=(n^{\alpha}\phi_{\xi_{n}})^{1/d}$,
$\xi_{n}$ will be chosen later. Set $u_{0}(x)=\sum_{n=n_{0}}^{\infty}u_{n}(x)$, where $n_{0}$ satisfies $\frac{1}{\beta_{n_{0}}}<\frac{1}{3}R$. For all $n\geq n_{0}$, $B_{1/\beta_{n}+2/3R}\subset B_{R}$. Note that
\begin{align*}
\|u_{0}\|_{L^{1}}\leq \sum_{n=1}^{\infty}\|u_{n}\|_{L^{1}}=\omega_{d}\sum_{n=1}^{\infty}n^{-\alpha}<\infty.
\end{align*}
Similar to the argument in the proof of Theorem \label{local}, we have
\begin{align*}
\int_{B_{R}}u(t,x)dx\geq \int_{B_{R}}\int_{0}^{t}S_{\alpha}(t-s)f(S_{\alpha}(s)u_{n})dsdx.
\end{align*}
Let $v_{0}=\psi\beta^{d}\chi_{1/\beta}$, $\psi$ is a positive constant to be chosen later. With $r=\frac{1}{\beta}$ and $\delta=\frac{R}{3}$,  by Lemma \ref{constant}, we get
\begin{align}\label{v}
S_{\alpha}(s)v_{0}\geq c(d,\alpha)\frac{\psi \beta^{d}}{(1+\beta s^{1/\alpha})^{d}}\chi_{1/\beta+t^{1/\alpha}}.
\end{align}
From (\ref{v}), it follows that for $|x|\leq \frac{1}{\beta}+t^{1/\alpha}$ and $s\leq t_{k}:=\min\{(\frac{R}{3})^{\alpha}, [(\frac{\psi}{\phi_{k}})^{1/d}-\frac{1}{\beta}]^{\alpha}\}$,
\begin{align}\label{svcd}
S_{\alpha}(s)v_{0}\geq c(d,\alpha)\phi_{k},
\end{align}
where $\psi$ satisfies $\psi\geq \frac{\phi_{k}}{\beta^{d}}$.

For $t\in (0,(\frac{R}{3})^{\alpha})$, since $f$ is non-decreasing, by (\ref{svcd}) and Corollary \ref{integrate},
\begin{align}\label{sum}
&\int_{B_{R}}\int_{0}^{t}S_{\alpha}(t-s)f(S_{\alpha}(s)v_{0})dsdx\nonumber\\
&\geq \sum_{k}\int_{B_{R}}\int_{t_{k+1}}^{t_{k}}S_{\alpha}(t-s)f(S_{\alpha}(s)v_{0})dsdx\nonumber\\
&=\sum_{k}\int_{t_{k+1}}^{t_{k}}\int_{B_{R}}S_{\alpha}(t-s)f(S_{\alpha}(s)v_{0})dxds\nonumber\\
&\geq \sum_{k}f(c(d,\alpha)\phi_{k})\int_{t_{k+1}}^{t_{k}}\int_{B_{R}}S_{\alpha}(t-s)\chi_{1/\beta+t^{1/\alpha}}dxds\nonumber\\
&\geq \mu(d,\alpha) \sum_{k}f(c(d,\alpha)\phi_{k})\int_{t_{k+1}}^{t_{k}}(\frac{1}{\beta}+s^{1/\alpha})^{d}ds\nonumber\\
&\geq \mu(d,\alpha) \sum_{k}f(c(d,\alpha)\phi_{k})\int_{t_{k+1}}^{t_{k}}s^{d/\alpha}ds,
\end{align}
where the sum in $k$ is taken over those for which
\begin{align}\label{those}
\frac{1}{\beta^{d}}\phi_{k}\leq \psi \leq (\frac{1}{\beta}+t^{1/\alpha})^{d}\phi_{k}.
\end{align}
For $k$ satisfying (\ref{those}) and $\frac{\alpha}{\beta^{d}}\phi_{k+1}\leq \psi$, we obtain
\begin{align}\label{compute}
\int_{t_{k+1}}^{t_{k}}s^{d/\alpha}ds&=\frac{\alpha}{d+\alpha}(t_{k}^{d/\alpha+1}-t_{k+1}^{d/\alpha+1})\nonumber\\
&=\frac{\alpha}{d+\alpha}\big\{[(\frac{\psi}{\phi_{k}})^{1/d}-\frac{1}{\beta}]^{d+\alpha}-[(\frac{\psi}{\phi_{k+1}})^{1/d}-\frac{1}{\beta}]^{d+\alpha}\big\}\nonumber\\
&=\frac{\alpha}{d+\alpha}(\frac{\psi}{\phi_{k}})^{1+\alpha/d}\big\{\big[1-(\frac{\phi_{k}}{\psi \beta^{d}})^{1/d}\big]^{d+\alpha}-\frac{\phi_{k}^{1+\alpha/d}}{\phi_{k+1}^{1+\alpha/d}}\big[1-(\frac{\phi_{k+1}}{\psi \beta^{d}})^{1/d}\big]^{d+\alpha}\big\}\nonumber\\
&\geq \frac{\alpha}{d+\alpha}(\frac{\psi}{\phi_{k}})^{1+\alpha/d}\big(1-\frac{\phi_{k}^{1+\alpha/d}}{\phi_{k+1}^{1+\alpha/d}}\big)\big[1-(\frac{\phi_{k+1}}{\psi \beta^{d}})^{1/d}\big]^{d+\alpha}.
\end{align}
Since $\phi_{k+1}\geq \tau\phi_{k}$, $\frac{2}{\beta^{d}}\phi_{k+1}\leq \psi$, by (\ref{compute}), we have
\begin{align}\label{certain}
\int_{t_{k+1}}^{t_{k}}s^{d/\alpha}ds\geq c(d,\alpha,\tau)\frac{\alpha}{d+\alpha}(\frac{\psi}{\phi_{k}})^{1+\alpha/d},
\end{align}
where $c(d,\alpha,\tau)$ is a positive constant.

By (\ref{sum}), (\ref{certain}),
\begin{align}\label{sufficient}
\int_{B_{R}}\int_{0}^{t}S_{\alpha}(t-s)f(S_{\alpha}(s)v_{0})ds
&\geq \mu(d,\alpha)c(d,\alpha,\tau)\sum_{k}f(c(d,\alpha)\phi_{k})(\frac{\psi}{\phi_{k}})^{1+\alpha/d}\nonumber\\
&=c\psi^{p}\sum_{k}f(s_{k})s_{k}^{-p},
\end{align}
where the sum in $k$ is taken over
\begin{align}\label{set}
\big\{k: \frac{\alpha}{\beta^{d}}\leq \frac{\psi}{\phi_{k+1}}\leq \frac{\psi}{\phi_{k}}\leq (\frac{1}{\beta}+t^{1/\alpha})^{d}\big\}.
\end{align}
For $t\in (0, (\frac{R}{3})^{\alpha})$, when $n$ is large enough such that $t^{1/\alpha}n^{\alpha/d}\geq 1$, let $\psi=n^{-\alpha}$, $\beta=\beta_{n}=(n^{\alpha}\phi_{\xi_{n}})^{1/d}$, then the set (\ref{set}) contains
\begin{align*}
\big\{k: 1\leq \phi_{k}, \mbox{and} \ \ \phi_{k+1} \leq\frac{1}{\alpha}\phi_{\xi_{n}}\big\}=\{k:k_{0}\leq k\leq k_{n}\},
\end{align*}
where $k_{0}$ is the smallest value of $k$ such that $\phi_{k}\geq 1$, $\xi_{n}$ is chosen to satisfy $\phi_{k_{n}+1}\leq \frac{1}{\alpha}\phi_{\xi_{n}}$.
Since $\sum_{k=1}^{\infty}s_{k}^{-p}f(s_{k})=\infty$, by (\ref{set}), $\phi_{k}=\frac{s_{k}}{c(d,\alpha)}$ and $\psi=n^{-\alpha}$, we can choose $k_{n}$ such that $n^{-\alpha p}\sum_{k=k_{0}}^{k_{n}}f(s_{k})s_{k}^{-p}\rightarrow \infty$ as $n\rightarrow \infty$.
\end{proof}

We can get an equivalent integral condition for non-existence of a local $L^{1}$ solution of (\ref{L1}).
\begin{lemma}\label{add} (See Lemma 4.2 in \cite{LRSV}.)
Assume that $f:[0,\infty)\rightarrow [0,\infty)$ is continuous and non-decreasing and $p>1$. The following conditions are equivalent.

(i) The exists a sequence $\{s_{k}\}_{k=1}^{\infty}$ such that $s_{k+1}\geq \tau s_{k}$, $\tau>1$ and
\begin{align*}
\sum_{k=1}^{\infty}s_{k}^{-p}f(s_{k})=\infty.
\end{align*}

(ii) $\int_{1}^{\infty}s^{-p}F(s)ds=\infty$, where $F(s)=\sup_{1\leq t\leq s}\frac{f(t)}{t}$.
\end{lemma}

\begin{remark}
In the Lemma \ref{add}, $f$ is supposed to be continuous. In fact, in Lemma 4.2 in \cite{LRSV}, there is no such requirement. But in the proof of  Lemma 4.2 in \cite{LRSV}, some sequence is assumed to make that the supreme in (ii) is attainable, this needs some additional conditions other than $f$ is non-decreasing.
\end{remark}

\subsection{An integral condition for local existence}
\begin{lemma}\label{super}
Let $u_{0}\geq 0$. If $f: [0,\infty)\rightarrow [0,\infty)$ is continuous and non-decreasing and there exists a function $v\in L^{1}((0,T)\times B_{R})$ such that
\begin{align}
S_{\alpha}(t)u_{0}+\int_{0}^{t}S_{\alpha}(t-s)f(v(s))ds\leq v(t), \ t\in [0,T],
\end{align}
then there exists a local integral solution $u$ of (\ref{sfh}), $u(t,x)\leq v(t,x)$, $x \in B_{R}$, $t\in [0,T]$.
\end{lemma}
\begin{proof}
By Lemma \ref{solution}, the conclusion holds.
\end{proof}

Let $X$, $Y$ be two Banach spaces. By $L(X,Y)$, we denote the space of all bounded linear operators from $X$ into $Y$. If $X=Y$, we simply rewrite $L(X,Y)$ as $L(X)$. A function $S: [0, \infty)\rightarrow L(X,Y)$ is strongly continuous if $t\mapsto S(t)x$ is continuous for all $x\in X$.
\begin{lemma}\label{convolution} (See Proposition 1.3.4 in \cite{ABHN}) Let $f \in L^{1}_{loc}([0,\infty), X)$ and let $S: [0,\infty)\rightarrow L(X,Y)$ be strongly continuous. Then the convolution
\begin{align*}
(S\ast f)(t)=\int_{0}^{t}S(t-s)f(s)ds
\end{align*}
exists (as a Bochner integral) and defines a continuous function $S\ast f: [0,\infty)\rightarrow Y$.
\end{lemma}

\begin{theorem}\label{L1existence}
If $f: [0,\infty)\rightarrow [0,\infty)$ is continuous and non-decreasing and
\begin{align}\label{integral condition}
\int_{1}^{\infty}s^{-(1+\alpha/d)}F(s)ds<\infty, \ where \ F(s)=\sup_{1\leq t\leq s}\frac{f(t)}{t},
\end{align}
then for every non-negative $u_{0}\in L^{1}(B_{R})$, there exists a constant $T>0$ such that (\ref{L1}) has a solution
\begin{align*}
u\in L^{\infty}_{loc}((0,T); L^{\infty}(B_{R}))\cap C([0,T]; L^{1}(B_{R})).
\end{align*}
\end{theorem}

\begin{proof}
First, we consider the case that $u_{0}=0$. We have
\begin{align}\label{ball}
S_{\alpha}(t)\chi_{B_{R}}&=\int_{B_{R}}p_{D}(t,x,y)\chi_{B_{R}}(y)dy\nonumber\\
&\leq \int_{R^{d}}p_{D}(t,x,y)\chi_{B_{R}}(y)dy\nonumber\\
&\leq \chi_{B_{R}}.
\end{align}
For $t\in [0,t_{0}]$, since $S_{\alpha}(t)$ is strongly continuous, there exists a positive constant $M_{t_{0}}$,
such that
\begin{align}\label{mt0}
\|S_{\alpha}(t)\|\leq M_{t_{0}}.
\end{align}
By ( \ref{ball}) and (\ref{mt0}), we obtain
\begin{align*}
S(t)u_{0}+\int_{0}^{t}S_{\alpha}(t-s)f(S_{\alpha}(s)\chi_{B_{R}})ds&\leq \int_{0}^{t}S_{\alpha}(t-s)(f(1)\chi_{B_{R}})ds\\
&\leq tM_{t_{0}}f(1)\chi_{B_{R}}\\
&\leq \chi_{B_{R}}
\end{align*}
for all $t$ small enough.

Second, we consider the case that $u_{0}\neq0$. Define $\tilde{f}(s)=f(s)$ for $s\in [0,1]$ and $\tilde{f}(s)=sF(s)$ for $s>1$. We see that $f(s)\leq \tilde{f}(s)$, $\tilde{f}(s)/s: [1,\infty)\rightarrow [0,\infty)$ is non-decreasing. Since any supersolution of the equation
\begin{align}\label{tilde}
u_{t}=\Delta^{\alpha/2}u+\tilde{f}(u)
\end{align}
is a supersolution of (\ref{L1}). We only need to show that (\ref{tilde}) has a supersolution.

The integral condition in (\ref{integral condition}) can be rewritten as
\begin{align*}
\int_{1}^{\infty}s^{-(2+\alpha/d)}\tilde{f}(s)ds<\infty,
\end{align*}
letting $s=\tau^{-d/\alpha}$, we have
\begin{align*}
\int_{0}^{1}\tau^{d/\alpha}\tilde{f}(\tau^{-d/\alpha})d\tau<\infty.
\end{align*}
For any $B>1$, set  $v(t)=BS(t)u_{0}+\chi_{B_{R}}$. Then
\begin{align}\label{B}
\mathcal{F}(v(t))&=S_{\alpha}(t)u_{0}+\int_{0}^{t}S_{\alpha}(t-s)\tilde{f}(v(s))ds\nonumber\\
&=S_{\alpha}(t)u_{0}+\int_{0}^{t}S_{\alpha}(t-s)\tilde{f}(BS_{\alpha}(s)u_{0}+\chi_{B_{R}})ds\nonumber\\
&=S_{\alpha}(t)u_{0}+\int_{0}^{t}S_{\alpha}(t-s)\frac{\tilde{f}(BS_{\alpha}(s)u_{0}+\chi_{B_{R}})}{BS_{\alpha}(s)u_{0}+\chi_{B_{R}}}(BS_{\alpha}(s)u_{0}+\chi_{B_{R}})ds\nonumber\\
&\leq S_{\alpha}(t)u_{0}+\int_{0}^{t}S_{\alpha}(t-s)\|\frac{\tilde{f}(BS_{\alpha}(s)u_{0}+\chi_{B_{R}})}{BS_{\alpha}(s)u_{0}+\chi_{B_{R}}}\|_{L^{\infty}}(BS_{\alpha}(s)u_{0}+\chi_{B_{R}})ds.
\end{align}
Since for all $t>0$, $S_{\alpha}(t)\chi_{B_{R}}\leq \chi_{B_{R}}$, by (\ref{B}),
\begin{align*}
\mathcal{F}(v(t))\leq S_{\alpha}(t)u_{0}+(BS_{\alpha}(t)u_{0}+\chi_{B_{R}})\int_{0}^{t}\|\frac{\tilde{f}(BS_{\alpha}(s)u_{0}+\chi_{B_{R}})}{BS_{\alpha}(s)u_{0}+\chi_{B_{R}}}\|_{L^{\infty}}ds.
\end{align*}
Since $\tilde{f}(s)/s$ is non-decreasing for $s\geq 1$, we get
\begin{align}\label{non-decreasing}
\mathcal{F}(v(t))\leq S_{\alpha}(t)u_{0}+(BS_{\alpha}(t)u_{0}+\chi_{B_{R}})\int_{0}^{t}\frac{\tilde{f}(\|BS_{\alpha}(s)u_{0}+\chi_{B_{R}}\|_{L^{\infty}})}{\|BS_{\alpha}(s)u_{0}+\chi_{B_{R}}\|_{L^{\infty}}}\|ds.
\end{align}
By the smoothing estimate (\ref{smooth effect}), we have
\begin{align}\label{estimate}
\|S_{\alpha}(t)u_{0}\|_{L^{\infty}}\leq Ct^{-d/\alpha}\|u_{0}\|_{L^{1}}.
\end{align}
Then
\begin{align}\label{bounded}
\|BS_{\alpha}(s)u_{0}+\chi_{B_{R}}\|_{L^{\infty}}\leq B\|S_{\alpha}(s)u_{0}\|_{L^{\infty}}+1\leq BCs^{-d/\alpha}\|u_{0}\|_{L^{1}}+1\leq 2BCs^{-d/\alpha}\|u_{0}\|_{L^{1}}
\end{align}
holds for small enough $s$.  \\
By (\ref{non-decreasing}) and (\ref{bounded}),  we have
\begin{align}\label{norm}
\mathcal{F}(v(t))\leq S_{\alpha}(t)u_{0}+(BS_{\alpha}(t)u_{0}+\chi_{B_{R}})\int_{0}^{t}\frac{\tilde{f}(2BCs^{-d/\alpha}\|u_{0}\|_{L^{1}})}{2BCs^{-d/\alpha}\|u_{0}\|_{L^{1}}}ds.
\end{align}
Note that $B>1$, letting $2BCs^{-d/\alpha}=\tau^{-d/\alpha}$, by (\ref{norm}), we obtain
\begin{align*}
\mathcal{F}(v(t))&\leq S_{\alpha}(t)u_{0}+\big(\int_{0}^{t(2BC\|u_{0}\|_{L^{1}})^{-\alpha/d}}\tau^{d/\alpha}\tilde{f}(\tau^{-d/\alpha})d\tau\big)(2BC\|u_{0}\|_{L^{1}})^{\alpha/d}(BS_{\alpha}(t)u_{0}+\chi_{B_{R}})\\
&\leq BS_{\alpha}(t)u_{0}+\chi_{B_{R}},
\end{align*}
for $t$ small enough. By Lemma \ref{super}, we obtain the local existence of a solution $u(t)$ which satisfies $u(t)\leq v(t)=BS_{\alpha}(t)u_{0}+\chi_{B_{R}}$ and $u\in L^{\infty}_{loc}((0,T); L^{\infty}(B_{R}))$. Since $f(u)\in L^{1}((0,T); L^{1}(B_{R}))$ and
\begin{align*}
u(t)=S_{\alpha}(t)u_{0}+\int_{0}^{t}S_{\alpha}(t-s)f(u(s))ds,
\end{align*}
Since $f(u)\in L^{1}((0,T); L^{1}(B_{R}))$, by Lemma \ref{convolution}, we have $u \in C([0,T]; L^{1}(B_{R}))$.
\end{proof}

\begin{corollary}
Let $f: [0,\infty)\rightarrow [0,\infty)$ is continuous and non-decreasing. Then (\ref{L1}) has the local $L^{1}$ existence property if and only if
\begin{align*}
\int_{1}^{\infty}s^{-(1+\alpha/d)}F(s)ds<\infty, \ where \ F(s)=\sup_{1\leq t\leq s}\frac{f(t)}{t}.
\end{align*}
\end{corollary}

\section{Results for the whole space $\mathbb{R}^{d}$}
Consider the equation (\ref{sfh}) on $\mathbb{R}^{d}$. Under the additional condition $\lim_{s\rightarrow 0}\frac{f(s)}{s}<\infty$, we have the following theorem.
\begin{theorem}
Let $\alpha\in (1,2]$ and let $f: [0,\infty)\rightarrow [0,\infty)$ be continuous and non-decreasing. Then \\
(i) for $q\in (1,\infty)$, (\ref{sfh}) has the local existence property in $L^{q}(\mathbb{R}^{d})$ if and only if
\begin{align*}
\lim_{s\rightarrow 0}\sup\frac{f(s)}{s}<\infty \ and \ \lim_{s\rightarrow \infty}s^{-(1+\alpha q)/d}f(s)<\infty;
\end{align*}
(ii)  (\ref{sfh}) has the local existence property in $L^{1}(\mathbb{R}^{d})$ if and only if
\begin{align*}
\lim_{s\rightarrow 0}\sup\frac{f(s)}{s}<\infty \ and \ \ \int_{1}^{\infty}s^{-(1+\alpha/d)}F(s)ds<\infty,
\end{align*}
where $F(s)=\sup_{1\leq t\leq s}\frac{f(t)}{t}$.
\end{theorem}
\begin{proof}
(i) We first show that if (\ref{sfh}) has the local existence property in $L^{q}(\mathbb{R}^{d})$, then the inequalities in (i) hold.
The non-existence results of Theorem \ref{local solution} still hold for the equations on $\mathbb{R}^{d}$, since proofs of them only need the lower bounds for the Dirichlet fractional heat kernel and  note that the inequality (\ref{Chen}) is similar to the inequality on the left side of (\ref{c1}). We will prove that if
$\lim_{s\rightarrow 0}\sup\frac{f(s)}{s}=\infty$, then there exists a non-negative $u_{0}\in L^{q}(\mathbb{R}^{d})$ such that the solution of (\ref{sfh}) is not bounded in $L^{q}(\mathbb{R}^{d})$. In fact, we can choose a sequence $s_{n}\rightarrow 0$ such that $s_{n}\leq n^{-\alpha}$ and $f(s_{n})\geq n^{2\alpha}s_{n}$. Set
\begin{align*}
u_{0}=\sum_{n=1}^{\infty}u_{n}, \ u_{n}=\frac{1}{\nu(d,\alpha)}s_{n}\chi_{n^{-\alpha q/d}s_{n}^{-q/d}}(x_{n}),
\end{align*}
where $\nu(d,\alpha)$ is the same constant as that in Corollary \ref{for which} and $x_{n}$ are chosen such that $B(x_{n}, n^{-\alpha q/d}s_{n}^{-q/d})$ are disjoint. Note that
\begin{align*}
\|u_{n}\|_{L^{q}}^{q}=\frac{1}{\nu(d,\alpha)^{q}}\omega_{d}n^{-\alpha q},
\end{align*}
then
\begin{align*}
\|u_{0}\|_{L^{q}}\leq \sum_{n=1}^{\infty}\|u_{n}\|_{L^{q}}=\frac{1}{\nu(d,\alpha)}\omega_{d}^{1/q}\sum_{n=1}^{\infty}n^{-\alpha}<\infty.
\end{align*}
Since $s_{n}\leq n^{-\alpha}$, then $n^{-\alpha q/d}s_{n}^{-q/d}\geq 1$. For all $0<s\leq 1$, by Corollary \ref{for which},
\begin{align*}
S_{\alpha}(s)u_{n}&=S_{\alpha}(s)\big(\frac{1}{\nu(d,\alpha)}s_{n}\chi_{n^{-\alpha q/d}s_{n}^{-q/d}}(x_{n})\big)\\
&=\frac{s_{n}}{\nu(d,\alpha)}S_{\alpha}(s)\chi_{n^{-\alpha q/d}s_{n}^{-q/d}}(x_{n})\\
&\geq s_{n}\chi_{n^{-\alpha q/d}s_{n}^{-q/d}+\sqrt{s}}(x_{n})\\
&\geq s_{n}\chi_{n^{-\alpha q/d}s_{n}^{-q/d}}(x_{n}).
\end{align*}
Thus
\begin{align}\label{un}
f(S_{\alpha}(s)u_{n})\geq f(s_{n})\chi_{n^{-\alpha q/d}s_{n}^{-q/d}}(x_{n}), \ \ 0<s\leq 1.
\end{align}
For all $0<t\leq1$, by (\ref{un}) and Corollary \ref{for which}, we obtain
\begin{align}\label{together}
S_{\alpha}(t-s)f(S_{\alpha}(s)u_{n})&\geq \nu(d,\alpha)f(s_{n})\chi_{n^{-\alpha q/d}s_{n}^{-q/d}+(t-s)^{1/\alpha}}(x_{n})\nonumber\\
&\geq \nu(d,\alpha)f(s_{n})\chi_{n^{-\alpha q/d}s_{n}^{-q/d}}(x_{n}).
\end{align}
Therefore,
\begin{align*}
u(t,x)&\geq \int_{0}^{t}S_{\alpha}(t-s)f(u(s))ds\\
&\geq \int_{0}^{t}S_{\alpha}(t-s)f(S_{\alpha}(s)u_{0})ds\\
&\geq \int_{0}^{t}S_{\alpha}(t-s)f(S_{\alpha}(s)u_{n})ds\\
&\geq \nu(d,\alpha)\int_{0}^{t}f(s_{n})\chi_{n^{-\alpha q/d}s_{n}^{-q/d}}(x_{n})ds\\
&\geq \nu(d,\alpha)tn^{2\alpha}s_{n}\chi_{n^{-\alpha q/d}s_{n}^{-q/d}}(x_{n}),
\end{align*}
for all $0<t\leq 1$. \\
Then
\begin{align*}
\|u(t)\|_{L^{q}}^{q}=\int_{\mathbb{R}^{d}}|u(t,x)|^{q}dx\geq \nu(d,\alpha)^{q}t^{q}\omega_{d}n^{\alpha q}\rightarrow \infty
\end{align*}
as $n\rightarrow \infty$.

Second, for $q\in (1,\infty)$, we will prove that if
\begin{align}\label{add}
\lim_{s\rightarrow 0}\sup\frac{f(s)}{s}<\infty \ and \ \lim_{s\rightarrow \infty}s^{-(1+\alpha q)/d}f(s)<\infty,
\end{align}
then (\ref{sfh}) has the local existence property in $L^{q}(\mathbb{R}^{d})$.  \\
For convenience, we let $p=1+\alpha q/d$. From (\ref{add}), it follows that
\begin{align*}
f(s)\leq C(s+s^{p}),
\end{align*}
where $C$ is a positive constant. When $f(s)=2Cs^{p}$,  (\ref{sfh}) has the local $L^{q}$ existence property for $t$ is sufficiently small, see \cite{ZY}. Given a non-negative $u_{0}\in L^{q}(\mathbb{R}^{d})$, let $u(t)$ be the local $L^{q}$ solution with $f(s)=2Cs^{p}$. Set $v(t)=e^{2Ct}u(t)$, we have
\begin{align}\label{supersolution}
v_{t}+(-\triangle)^{\alpha/2}v&=2Ce^{2Ct}u+e^{2Ct}(u_{t}+(-\triangle)^{\alpha/2}u)\nonumber\\
&=2Ce^{2Ct}u+e^{2Ct}2Cu^{p}\nonumber\\
&\geq 2C(v+(e^{2Ct})^{1-p}v^{p})\nonumber\\
&\geq 2Cv+v^{p}>C(v+v^{p}),
\end{align}
where $(e^{2Ct})^{1-p}\geq 1/2$.

From $(e^{2Ct})^{1-p}\geq 1/2$, we have $0<t\leq \frac{\ln2}{2C(p-1)}$. By (\ref{supersolution}), it is easy to see that $v$ is a supersolution for $0<t\leq \frac{\ln2}{2C(p-1)}$. Then from Lemma \ref{super}, it follows that there exists a solution of (\ref{sfh}) bounded in $L^{q}(\mathbb{R}^{d})$.

(ii) The non-existence results of Theorem \ref{assume} still hold for the equations on $\mathbb{R}^{d}$, by the same argument as that in the proof of (i). By the simliar proof as that in (i), we can prove that if
$\lim_{s\rightarrow 0}\sup\frac{f(s)}{s}=\infty$, then there exists a non-negative $u_{0}\in L^{1}(\mathbb{R}^{d})$ such that the solution of (\ref{sfh}) is not bounded in $L^{1}(\mathbb{R}^{d})$. Next we show that if
\begin{align*}
\lim_{s\rightarrow 0}\sup\frac{f(s)}{s}<\infty \ \mbox{and} \ \ \int_{1}^{\infty}s^{-(1+\alpha/d)}F(s)ds<\infty,
\end{align*}
where $F(s)=\sup_{1\leq t\leq s}\frac{f(t)}{t}$,
then  (\ref{sfh}) has the local existence property in $L^{1}(\mathbb{R}^{d})$. In fact, by $\lim_{s\rightarrow 0}\sup\frac{f(s)}{s}<\infty$ and the continuity of $f$, we have $f(0)=0$.  If $u_{0}=0$, then $u(t)\equiv 0$ is a solution of (\ref{sfh}). If $u_{0}\neq 0$, for $x\in\mathbb{R}^{d}$,
\begin{align*}
(S_{\alpha}(t)u)(x)=\int_{\mathbb{R}^{d}}p(t,x,y)u(y)dy,
\end{align*}
we get
$S_{\alpha}(t)1=1$ for all $t>0$, and $S_{\alpha}(t)u_{0}>0$ for $u_{0}>0$.  Then by similar arguments to that in the proof of Theorem \ref{L1existence}, we obtain the local $L^{1}$ existence property.
\end{proof}
\section{Acknowledgements}

Kexue Li is supported by National Natural Science Foundation of China under the contract
No.11571269, China Postdoctoral Science Foundation Funded Project under the contract No.2015M572539 and Shaanxi Province Postdoctoral Science Foundation Funded Project. He would like to thank China Scholarship Council that has provided a
scholarship for his research work in the United States.

\end{document}